\documentclass[a4paper,12pt,leqno]{amsart}

\usepackage{amssymb,hyperref}
\usepackage[latin1]{inputenc}
\theoremstyle{plain}
\newtheorem{theo}{Theorem}
\theoremstyle{definition}
\newtheorem{quest}{Question}
\theoremstyle{remark}

\numberwithin{equation}{section}

\newcommand{\thismonth}{\ifcase\month\or
  January\or February\or March\or April\or May\or June\or July\or
  August\or September\or October\or November\or December\fi
  \space\number\year}
\DeclareMathAlphabet{\mathrmsl}{OT1}{cmr}{m}{sl}

\newcommand{\oper}[3][n]{\newcommand{#2}{\mathop{\mathrm{#3}}%
\ifx n#1\nolimits\else\limits\fi} }
\newcommand{\rsoper}[3][n]{\newcommand{#2}{\mathop{\mathrmsl{#3}}%
\ifx n#1\nolimits\else\limits\fi} }
\newcommand{\TM}{{\mathbb T}}
\newcommand{\PM}{{\mathbf P}}
\newcommand{\RM}{{\mathbb R}}
\newcommand{\CM}{{\mathbb C}}

\newcommand{\vol}{\operatorname{vol}}
\newcommand{\Scal}{\operatorname{Scal}}

\newcommand{\End}{\operatorname{End}}

\renewcommand{\geq}{\geqslant}
\renewcommand{\leq}{\leqslant}
\newcommand{\cerc}{{\mathbb S}}

\newcommand{\proofof}[1]{\end{#1}\begin{proof}}

\newcounter{mnotecount}[section]
\renewcommand{\themnotecount}{\thesection.\arabic{mnotecount}}
\newcommand{\mnote}[1]%{}
{\protect{\stepcounter{mnotecount}}$^{\mbox{\footnotesize  $%\!\!\!\!\!\!\,
      \bullet$\themnotecount}}$ \marginpar{\raggedright\tiny\em
    $\!\!\!\!\!\!\,\bullet$\themnotecount: #1} }

\begin{document}

\title[Vafa--Witten bound on the sphere]{Extremality for the 
Vafa--Witten bound\\ on the sphere}
\author{Marc Herzlich}
\address{Institut de Math\'ematiques et Mod\'elisation de Montpellier\\ 
UMR 5149 du CNRS\\ Université Mont­pellier~II\\ France}
\email{herzlich@math.univ-montp2.fr}

\begin{abstract}
We prove that the round metric on the sphere has the largest first eigenvalue
of the Dirac operator among all metrics that are larger than it. As a corollary,
this gives an alternative proof of an extremality result for scalar curvature
due to M. Llarull. 
\end{abstract}
\thanks{The author is a member of the {\sc edge} Research Training 
Network {\sc hprn-ct-2000-00101} of the European Union  and is supported 
in part by an {\sc aci} program of the French Ministry of Research.}

\maketitle

\section{Introduction}

\smallskip

C. Vafa and E. Witten have shown in \cite{vafawitten} that there exists 
a common upper bound for the smallest eigenvalue of all {\sl twisted} Dirac 
operators on a given Riemannian manifold. The common upper bound they found 
depends strongly on the choice of the metric on the base. Using an analogous 
method, H. Baum \cite{baum-upper} exhibited later an explicit upper bound for 
the first eigenvalue of the (untwisted) Dirac operator on an even-dimensional 
Riemannian manifold that can be sent on a sphere by a map of (high) non-zero 
degree. This bound depends on the Lipschitz norm of the map from the manifold 
to the round sphere.

The goal of this short note is to elaborate further on H. Baum's results 
\cite{baum-upper}. Indeed, we prove below that some {\sl optimal} upper 
bound can be obtained for the bottom of the spectrum of the (untwisted) Dirac 
operator of a large class of metrics on the sphere. This extends \cite{baum-upper}
in two different ways, first by providing an optimal bound, and secondly by 
extending it to odd dimensions.

\medskip

\begin{theo}
Let $g$ be any Riemannian metric and $b$ be the round metric of constant
curvature $1$ on the sphere $\cerc^{n}$. If $g\geq b$ pointwise, then 
there is an eigenvalue of the Dirac operator of $g$ in $[-\frac{n}{2},\frac{n}{2}]$.
Moreover, if no eigenvalue lies in the open interval 
$]-\frac{n}{2},\frac{n}{2}[$, then $g$ is isometric to the round metric $b$.
\end{theo}

\medskip

As an interesting corollary, one gets an alternative proof of M. Llarull's 
extremality result for scalar curvature on the sphere \cite{llarull-sphere}.
Indeed, the classical first eigenvalue estimate for the Dirac 
operator due to Th. Friedrich \cite{friedrich-dim3} states that the smallest 
eigenvalue  $\lambda_1(g)$ (in absolute value) of the Dirac operator
satisfies:
\[ | \lambda_1(g)|^2  \geq \frac{n}{4(n-1)}\inf \Scal_g .\]
Hence one concludes:

\medskip

\begin{theo}[Llarull \cite{llarull-sphere}]
Let $g$ be any metric and $b$ be the round metric of constant curvature
$1$ on the sphere $\cerc^{n}$. If $g\geq b$ pointwise, then
\[ \inf \Scal_g  \leq  n(n-1)\ , \]
the inequality being strict if $g$ is distinct from $b$.
\end{theo}

\medskip

The proof of Theorem 1 is close to the original approach of Vafa and Witten, 
see also M. Atiyah \cite{atiyah-vw} and H. Baum \cite{baum-upper}. Since we
look for optimal bounds, we have however to pay a special attention to 
all the estimates involved. As a result, this implies that two different proofs 
are required in 
the even- and odd-dimensional cases (see \cite{davaux-mathann} for another 
occurrence of this problem in a related context). The first one is obtained 
from an index argument, whereas the second one needs a slightly more delicate 
proof using spectral flow considerations.

\medskip

{\flushleft\it Acknowledgements}. It was Maung Min-Oo who pointed out
during a lecture in July 2001 \cite{minoo-castello} that Vafa--Witten was
optimal on a round sphere and that a strong link might exist with Llarull 
extremality. This idea appears (in a somewhat more hidden way) on pages
84--85 in M. Gromov's long and rich paper \cite{gromov-Kaire}. 
The author is then happy to acknowledge the origin of the 
ideas behind this note. He thanks Maung Min-Oo again, Andrei Moroianu and 
Uwe Semmelmann for useful comments, H\'el\`ene Davaux for her careful reading
of a first draft of this paper, and the referee for his remarks.

\bigskip

\section{Background material}

\smallskip

For any metric $g$ on $\cerc^{n}$, we let $\Sigma_g$ be its spin bundle. Using
the idea introduced by J. P. Bourguignon and P. Gauduchon \cite{jpbpg}, one may 
identify the spin bundle $\Sigma_g$ with the spin bundle $\Sigma_b$ through
a lift of the principal oriented orthonormal frame bundles isomorphism
\[ H : P_{SO}\cerc^{n}(b) \longrightarrow P_{SO}\cerc^{n}(g) \]
induced by the unique symmetric positive definite map $H$ such that
\[ g(H\cdot,H\cdot) = b(\cdot,\cdot) .\]
If $g\geq b$, this implies that $H$ has operator norm relative to $b$ 
bounded from above by $1$. 
This allows to transfer the Clifford action relative to $g$ on the bundle
$\Sigma_b$. Denoting by $\ell_b$ (resp. $\ell_g$)~: 
$T\cerc^n \rightarrow \End(\Sigma_b)$ 
the Clifford actions for $b$ and $g$, they are related as:
\[ \ell_b(\cdot) = \ell_g(H\cdot) \]
The Levi-Civita connection for any metric $g$ can be transferred in the same way
on $\Sigma_b$ as a metric connection (but with torsion). In what follows, we 
will always assume that all spinor bundles on the sphere have been identified
to those relative to $b$, as above.

\medskip

The following facts are extracted from classical texts on spinor geometry, see
also \cite{baum-killing,baum-upper,gilkey-book} which are useful references.

{\sl In even dimensions $n=2m$}, the spin bundle $\Sigma_g$ splits into two 
{\sl half-spin} bundles $\Sigma_g^+$ and $\Sigma_g^-$; on the sphere $\cerc^{2m}$, 
the whole spin bundle, or equivalently 
the sum of both half-spin bundles, is a trivial bundle obtained by restricting
the spinor bundle of $\RM^{2m+1}$ to the sphere; namely:
\[  \TM : = 
\Sigma_b = \Sigma^+_b \oplus \Sigma_b^- = \cerc^{2m} \times \CM^{2^m}.\]
Moreover, projections on each factor $\Sigma_b^{\pm}$ can be explicitly described 
at a point $x$ on the sphere as:
\[ \Pi_{\pm}(x) = \frac{1}{2} \left( 1 \pm i \ell_0(x) \right) \]
where $\ell_0$ is the algebraic Clifford action of $\CM l_{2m+1}$ on
$\CM^{2^m}$ (as a result, none of the half-spin bundles is trivial).
Using this expresion, it is easy to relate the trivial connection $\nabla^0$
(of flat space) acting on $\TM$ with the Levi-Civita connection of $b$ acting
on $\Sigma_b^{\pm}$. If $\psi$ 
is any spinor field in $\Sigma_b^{\pm}$, then:
\[ \nabla^{\pm}_X\psi = \Pi_{\pm}\nabla^0_X\Pi_{\pm} \psi = \nabla^0_X\psi
- \frac{1}{2}\, \ell_0(x) \ell_0(X) \psi  \]
and the full Levi-Civita connection for $b$ on $\Sigma_b = \Sigma_b^+ \oplus
\Sigma_b^-$ is $\nabla^b = \nabla^+ \oplus \nabla^-$ (note that actions $\ell_0$ 
on vectors tangent to $\cerc^n$ and $\ell_b$ may be identified). 

Moreover, the bundle $\TM\otimes\TM = \Sigma_b\otimes\Sigma_b$ is the (complex) 
differential form bundle $\Lambda^{\bullet}\cerc^n\otimes\CM$ of the sphere.
On this bundle, two different Clifford actions exist. The first is the usual
one for twisted spinor bundles, with the Clifford algebra acting on the left 
factor of the tensor product; we will continue to denote this one below by 
$\ell$, as it is defined on a decomposed element $\sigma\otimes\tau$ of 
$\Sigma_b\otimes\Sigma_b$ by:
\[ \ell_b(X)\left(\sigma\otimes\tau\right) = (\ell_b(X)\sigma)\otimes\tau . \]
But there is another, right-handed, one that we shall denote by $r$ and
which is defined by:
\[ r_b(X)\left(\sigma\otimes\tau\right) = \sigma\otimes (\ell_b(X)\tau). \]
Both Clifford actions can be explicitely described when one identifies
$\Sigma_b\otimes\Sigma_b$ with $\Lambda^{\bullet}\cerc^{2m}$ as follows: 
for any $p$-form $\omega$, and any $1$-form $\alpha$
\[ \ell_b(\alpha^{\sharp})\omega = \alpha\wedge\omega - 
i_{\alpha^{\sharp}}\omega, \ \ \ 
r_b(\alpha^{\sharp})\omega = (-1)^p \left( \alpha\wedge\omega 
+ i_{\alpha^{\sharp}}\omega\right) \]
(musical isomorphism $\sharp$ referring, as always, to $b$). For more information
on these points, the reader is referred to the book \cite{lawson-michelsohn}
(see also the paper \cite{dirac-pert} where a concise account of these facts is 
given). The reader should also be careful about tensor product connections, as they 
might not coincide: for instance, $\nabla^b\otimes 1 + 1\otimes\nabla^b$ is the 
Levi-Civita connection on differential forms on the sphere whereas 
$\nabla^0\otimes 1 + 1\otimes\nabla^0$ is the trivial connection induced 
from $\RM^{2m+1}$.

\medskip

{\sl In odd dimensions $n=2m-1$}, the Clifford algebra is isomorphic to the 
sum of two copies of the matrix algebra
$\textrm{End}(\CM^{2^{m-1}})$. Hence it has two inequivalent 
representations which lead to two different Clifford bundles $\Sigma_b$ and
$\Sigma_b'$ on the sphere $\cerc^{2m-1}$. They are equivalent as bundles with 
structure group $\textrm{Spin}(2m-1)$ but the Clifford actions differ.
Their sum is again a trivial bundle:  
\[ \TM :=\Sigma_b \oplus\Sigma_b' = \cerc^{2m-1} \times \CM^{2^m}\ ; \]
as above, it is the restriction to the sphere of the full spinor bundle 
of flat space $\RM^{2m}$. This bundle itself splits in two trivial
subbundles of rank $2^{m-1}$
\[ \TM = \TM^+ \oplus \TM^- \]
obtained by taking the $(\pm i)$-eigenspaces of the volume form of $\RM^{2m}$.
Here we have to take care that these
bundles are different from the bundles $\Sigma_b$ and $\Sigma_b'$ on the sphere 
already alluded to a few lines above: for instance, $\TM^+$ and $\TM^-$
are exchanged when Clifford multiplied by vectors, whereas $\Sigma_b$ and 
$\Sigma_b'$ are preserved by the Clifford action. Notice that, contrarily
to what happens if $n=2m$, the bundles $\TM^{\pm}$ are preserved by the 
flat connection $\nabla^0$. As final remark, we recall that $\TM\otimes\TM$ is 
again isomorphic to the full differential form bundle of $\RM^{2m}$.
 
\medskip

From now on and to avoid confusion, we will use notations $\TM$, $\Sigma_b$ and
$\Sigma_g$ in the following sense : they will always denote the corresponding bundle
relative to the metric $b$, but endowed either with the trivial 
connection $\nabla^0$, or the Levi-Civita connections $\nabla^b$ or $\nabla^{\pm}$ 
of $b$ on spinors, or the Levi-Civita connection $\nabla^g$ of $g$ on spinors.

\bigskip

\section{The proof: even-dimensional case}

\smallskip

We consider the bundle $S = \Sigma_g\otimes\Sigma_b^+$; as spin bundles
for different metrics always are identified, this means that we take the tensor 
product of the whole spin bundle $\Sigma_b$ with the half-spin $\Sigma^+_b$,
but endowed with the {\sl non-trivial} tensor product connection $\nabla^g\otimes
1 + 1\otimes\nabla^+$, where $\nabla^g$ is the Levi-Civita connection of the
metric $g$ and $\nabla^+$ has been defined in the previous section. 

We now apply Atiyah-Singer index theorem to the twisted Dirac operator
\begin{equation} \mathcal{D}^{g+} : \Sigma^+_g \otimes
\Sigma^+_b \longrightarrow \Sigma^-_g \otimes \Sigma^+_b .
\end{equation}
One of the most important consequences of index theory is the topological 
invariance of the index. The index of the following model (`round') Dirac 
operator:
\[  \mathcal{D}^{b+} : 
\Sigma_b^+\otimes\Sigma_b^+ \longmapsto \Sigma_b^-\otimes\Sigma_b^+ \] 
equals $1$ (it is the operator whose index equals
one half of the Euler number in dimension $4k+2$ and one half of the sum of 
the Euler number and the signature in dimension $4k$ ; hence this value
on the even-dimensional sphere, see for instance \cite[p.~95--96]{gilkey-book}
for an explicit derivation), so that
the index of the Dirac operator $\mathcal{D}^{g+}$ built with the connection 
$\nabla^g\otimes 1 + 1\otimes\nabla^+$ on the twisted spin bundle 
$S = \Sigma_g\otimes\Sigma_b^+$ is also $1$. Its kernel is then non-zero.

Let us now consider the tensor product bundle $\Sigma_g\otimes \TM$, endowed
this time with the connection $\nabla^g\otimes 1 + 1\otimes\nabla^0$. As the pair 
($\TM,\nabla^0)$ is a \emph{trivial} flat bundle on the sphere, the spectrum of 
the twisted Dirac operator attached to this connection is the same as the spectrum 
of the Dirac operator on $\Sigma_g$, but with each eigenvalue repeated $2^m$ times 
its multiplicity. Applying standard perturbation theory, one gets the first part of 
the statement of Theorem 1 if one is able to bound the difference 
\[ \bar{L} := \mathcal{D}^{g0} - 
\left( \mathcal{D}^{g+} \oplus \mathcal{D}^{g-}\right) : 
\Sigma_g\otimes\Sigma_b \longrightarrow \Sigma_g\otimes\Sigma_b. \]
This is easily computed on a decomposed section $\sigma\otimes\tau$
of $\Sigma_g\otimes\Sigma_b$:
if $\{e_i\}$ is any $g$-orthonormal basis on the sphere,
\begin{equation*}\begin{split} 
\bar{L} (\sigma\otimes\tau) & 
= \sum_{i=1}^n \, \ell_g(e_i) \left( \nabla^g_{e_i}\sigma\otimes\tau + 
   \sigma\otimes\nabla^0_{e_i}\tau - \nabla^g_{e_i}\sigma\otimes\tau
 - \sigma\otimes\nabla^+_{e_i}\tau \right) \\
& = \frac{1}{2} \, \sum_i \, \ell_g(e_i)\sigma \otimes \ell_0(x)\ell_0(e_i)\tau
\end{split}\end{equation*}
This shows that 
\begin{equation*} 
\bar{L} = \frac{1}{2}\, r_0(x) \, L,
\end{equation*}
where $L$ is the (pointwise) linear map of vector bundles given on a spinor 
field $\psi$ in $\Sigma_g\otimes\Sigma_b$ by 
\[ L \psi =  \sum_i \, \ell_g (e_i) r_b(e_i) \psi .\]
We now choose a $b$-orthonormal basis $\{\varepsilon_i\}$ that diagonalizes $H$.
Eigenvalues $\mu_i$ live in $]0,1]$ and $\{e_i = \mu_i \varepsilon_i\}$
is a $g$-orthonormal basis, and moreover,
\[ L\psi = \sum_i \, \mu_i^{-1}\, \ell_b(e_i) r_b(e_i) \psi 
= \sum_i \, \mu_i\, \ell_b(\varepsilon_i) r_b(\varepsilon_i) \psi  .\]
As Clifford multiplication by a $b$-unit vector has $b$-operator norm equal to 
$1$, one immediately obtains that 
\[ |L\psi|_b \ \leq \ \left(\sum \mu_i\right) |\psi|_b 
\ \leq \ n\, |\psi|_b \ , \]
with equality if and only if $\mu_i=1$ for all $i$. 
This concludes the proof of Theorem 1 in the even-dimensional case. \qed

\smallskip

{\flushleft\sl Remark}. Using the explicit expressions of $\ell_b$ and $r_b$ 
on differential forms, one gets
\begin{equation*}
\ell_b(\varepsilon_k)r_b(\varepsilon_k)\, 
(\varepsilon_{i_1}\wedge\cdots\wedge\varepsilon_{i_p}) = \begin{cases}
(-1)^p\ \varepsilon_{i_1}\wedge\cdots\wedge\varepsilon_{i_p} & \textrm{ if } 
k \in \{i_1,\ldots,i_p\}, \\ 
(-1)^{p+1}\ \varepsilon_{i_1}\wedge\cdots\wedge\varepsilon_{i_p} & \textrm{ if 
not}, 
\end{cases}
\end{equation*}
and this of course leads to an explicit form of the previous estimate.

\bigskip

\section{The proof: odd-dimensional case}

\smallskip

Our treatment of the odd-dimensional case is much inspired by 
\cite{atiyah-vw,davaux-mathann}. Starting from the Dirac operator
on the spin bundle $\Sigma_g$, our goal is to find a trivial twisting
bundle endowed with a non-flat connection whose associated twisted
Dirac operator has a non-zero kernel. As index considerations are
useless in odd dimensions, we have to rely on a spectral flow argument.

\medskip

Let the dimension be $n=2m-1$. Our choice of bundle on $\cerc^n$ is 
\[ \Sigma_g\otimes\TM^+ = \Sigma\otimes \CM^{2^{m-1}}\] 
(notations as in section 2), so that, as above, the 
Dirac operator $\mathcal{D}^{g0}$ associated to $\nabla^g\otimes 1 + 
1\otimes\nabla^0$ on $\Sigma_g\otimes\TM^+$ has the same spectrum as 
$\mathcal{D}^g$ on $\Sigma_g$ (up to multiplicity). Note that this
is the only possible choice of connection at this stage since $\TM^{\pm}$
is not preserved by the Levi-Civita connection of $(\cerc^n,b)$.

We now fix $e$ unit in $\RM^{2m}$
and we define for each $x$ in $\cerc^{2m-1}$ an endomorphism $u(x)$ of
the fiber of $\TM^+$ at $x$ as follows: we start with $ (e\cdot x\cdot)$
which is an element $\CM l_{2m}$ that preserves both factors $\TM^{\pm}$.
 Hence we let:
\[ u(x) = \textrm{ the projection of } (e\cdot x\cdot)\ 
\textrm{ in } \textrm{End}(\TM^+) \ \ ( = \textrm{End}(\CM^{2^{m-1}}) ) .\]
As each $x$ is skew-hermitian and of square $-1$, $u(x)$ is unitary
for each $x$
and one can define the unitarily (gauge) equivalent connection 
\[ \nabla^u =  u^{-1}\circ \nabla^0 \circ u  \]
on $\TM^+$. The path of connections
\[ t \in [0,1] \longmapsto \nabla^g\otimes 1 + 1\otimes\nabla^t \ \textrm{ with } 
\ \nabla^t = (1-t)\nabla^0 + t\nabla^u \ . \]
gives rise to a path of Dirac operators $\mathcal{D}^t$ on $\Sigma_g\otimes\TM^+$
with unitarily equivalent operators at $\{ t=0\}$ and $\{ t=1\}$. 
The spectral flow of this family can be computed using the index theorem  
on $\cerc^{2m-1}\times\cerc^1$, applied to the Dirac operator acting on the 
spinor bundle $\Sigma_b$ twisted by the bundle obtained by identifying $\TM^+$ 
at $\{ t=0\}$ and $\{ t=1\}$ through $u$.

This spectral flow is non-zero: 
as a matter of fact, its value is 
\[ sf \ = \ \int_{\cerc^{2m-1}\times\cerc^{1}} \hat{A}(\cerc^{2m-1}\times\cerc^{1})
\,\textrm{ch} (\widetilde{\TM}^+) \ \ \ ,\]
where $\widetilde{\TM}^+$ is the bundle obtained on $\cerc^{2m-1}\times\cerc^{1}$
from $\TM^+\times [0,1]$ by identifying through $u$ at $\{ t=0\}$ and $\{ t=1\}$.
Moreover,
\[ \hat{A}(\cerc^{2m-1}\times\cerc^{1}) = 1, \]
and $\textrm{ch} (\widetilde{\TM}^+)$ can be computed with the curvature $2$-form
\[ \Omega = dt\wedge u^{-1}du - (t-t^2) u^{-1}du\wedge u^{-1} du \]
obtained by writing the connection $\nabla^t = \nabla^0 + \omega^t$, where 
$\omega^t$ is the $1$-form with values in $\End(\TM^+)$ equal to 
$ t\, u^{-1}(du \cdot)$; hence
\begin{equation*}\begin{split} \textrm{ch}(\Omega) 
& \ = \ \left(\frac{i}{2\pi}\right)^m \textrm{Tr}(\Omega^m) \\
& \ = \ C\,(t-t^2)^{m-1} dt\wedge \textrm{Tr}(u^{-1}du\wedge\dots\wedge u^{-1}du)
\end{split}\end{equation*}
with $C$ a non-zero constant. This computation shows that 
\[ sf \ = \ \int_{\cerc^{2m-1}\times\cerc^{1}} \textrm{ch} (\Omega) 
\ = \ C\,\int_{[0,1]} (t-t^2)^{m-1}dt \int_{\cerc^{2m-1}} 
\textrm{Tr}(x\cdot e_1\cdots e_{2m-1})\, d\!\vol_{\cerc^{2m-1}} \]
and this last result is non-zero by the very definition of $\TM^{\pm}$ as the
$(\pm i)$-eigenspaces of the volume form of $\RM^{2m}$ (see 
\cite{davaux-these,davaux-mathann} for instance for more detailed computations
in an analogous case).

As a result, there exists a value of $t$ in $[0,1]$ for which the kernel of 
$\mathcal{D}^t$
is non-zero. Suppose $t\leq 1/2$, then, as in the previous section, what is needed 
to conclude is an estimate of the difference
\[ \bar{L}_t \ = \ \mathcal{D}^t - \mathcal{D}^{g0} \]
(in case $t>1/2$, on should use $u^{-1}\circ\mathcal{D}^{g0}\circ u$ rather than 
$\mathcal{D}^{g0}$ but this is harmless since both operators have the same 
spectrum). An easy computation
shows that this difference of operators acts as the following linear map on 
$\Sigma_g\otimes\TM$ seen as a sub-bundle of the whole differential form 
bundle $\TM\otimes\TM$ of $\RM^{2m}$:
\[ \bar{L}_t =  r_b(t\, x)\,\circ\left( \sum_{i=1}^{2m-1} \mu_i\, 
\ell_b(\varepsilon_i)\, r_b(\varepsilon_i)\,\right) \]
(same notations as in the previous section).
The end of the proof is then entirely analogous to that in even dimensions.
\qed

\bigskip

\section{Final comments}

\smallskip
 
It is known that scalar curvature extremality extends to a large family of 
Riemannian manifolds; among these, one can find complex projective spaces, 
K\"ahler manifolds of positive Ricci curvature, non-negatively curved locally 
symmetric spaces (see the foundational \cite[\S 6]{gromov-Kaire} and 
\cite{goette-semm-spinc,goette-semm-symm,kramer-fibrations,llarull-sphere,minoo-symmetric}). 
In view of Theorem 2, it seems then natural to ask the following obvious questions:

\medskip

\begin{quest} 
Can one find other examples of eigenvalue-extremal metrics~? 
\end{quest}

\smallskip

\begin{quest} In particular,
if $g$ is a metric on $\CM\PM^m$ (with $m$ odd, so that it is spin) that is larger 
than the Fubini-Study metric $f$, does one have $|\lambda_1(g)| \leq 
|\lambda_1(f)|$ with equality if and only if $g$ is isometric to $f$~?
\end{quest}

\bigskip

\bibliographystyle{amsplain}

\providecommand{\bysame}{\leavevmode\hbox to3em{\hrulefill}\thinspace}

\medskip

\end{document}